\documentclass[leqno,twoside,11pt]{amsart}

\usepackage{latexsym,esint,url}
\usepackage{graphicx}

\theoremstyle{plain}

\theoremstyle{definition}

\begin{document}

\title[Game Theoretical Methods in Nonlinear PDEs]
{Game Theoretical Methods in PDEs}
\author{Marta Lewicka and Juan J. Manfredi}
\address{Marta Lewicka and Juan J. Manfredi, University of Pittsburgh, Department of Mathematics, 
139 University Place, Pittsburgh, PA 15260}
\email{lewicka@pitt.edu, manfredi@pitt.edu}
%\keywords{Tug-of-War games, viscosity
%  solutions, p-laplacian, mean value property}
%\date{\today}

\maketitle

Nonlinear PDEs, mean value properties,  and stochastic
differential games are intrinsically connected. We will describe how
the solutions to certain PDEs (of $p$-Laplacian type) can be interpreted as
limits of values of a specific Tug-of-War game, when the step-size $\epsilon$ determining the allowed
length of move of a token, decreases to $0$. 
This approach  originated in \cite{PSSW} and \cite{PS};  for the case
of deterministic games see the review SIAM News article \cite{K} and
\cite{KS, KS2}. 
% This article is based on  the minitutorial given by
% Marta Lewicka at the 2013 SIAM Conference on Analysis of Partial
% Differential Equations in Lake Buena Vista, FL  last December
% \cite{L}, where the reader can find many more details. 

\par\bigskip

\noindent {\bf How the linear elliptic equations arise in
  probability.}

Let us begin with the case governed by the discrete Brownian motion. Consider an open bounded set
$\Omega\subset\mathbb{R}^N$  and a non-empty portion of its boundary
$\Gamma_1\subset\partial\Omega$.  Place a token at a point
$x_0\in\Omega$ and assume that at each step of the process,   it is moved with
equal probabilities $\frac{1}{2N}$, to one of the $2N$ symmetric positions
$x_0\pm\epsilon e_i$, $i:1\ldots N$. Denote by  $u_\epsilon
(x_0)$ the probability that the first time the token exists $\Omega$,
it exits across $\Gamma_1$.  By applying conditional probabilities, it is
clear that $u_\epsilon$ satisfy the {\em mean value property}:
\begin{equation}\label{mean1}
\frac{1}{2N} \sum_{i=1}^N \Big( u_\epsilon(x+ \epsilon e_i) +
  u_\epsilon(x-\epsilon e_i) \Big) = u_\epsilon(x)
\end{equation}
Further, it follows that as $\epsilon\to 0$, the functions
$u_\epsilon$ converge uniformly in $\Omega$ to a
continuous $u\in\mathcal{C}(\Omega)$ which is a {\em
  viscosity solution} to the problem $\Delta u = 0$ in $\Omega$,
$u=\chi_{\Gamma_1}$ on $\partial\Omega$. 

More precisely,
this means that: (i) for each $x_0\in\Omega$ and each smooth test function $\phi$ satisfying
$u(x) -\phi(x) > u(x_0)-\phi(x_0) = 0$ for all $x\neq x_0$ in a small
neighbourhood of $x_0$, one has: $\Delta\phi(x_0)\leq
0$, (ii)  the same condition holds if we replace $u$ by
$-u$. It is well known that the viscosity solutions to $\Delta u = 0$
coincide with the classical solutions. An avantage of working with
the above, seemingly, more complex notion, is that the limiting properties of
$u_\epsilon$ follow quite naturally from the mean value property
(\ref{mean1}). Namely,  replacing the increments $u_\epsilon(x\pm\epsilon
e_i) - u_\epsilon(x)$ in the discontinuous $u_\epsilon$ in
(\ref{mean1}), by the same increments in the smooth $\phi$, applying
Taylor's expansion and taking into account the assumed sign of $u-\phi$, yields the sign of
$\Delta\phi$, wheras the first derivatives cancel out due to the symmetry
in (\ref{mean1}).

Heuristically,  this  can be seen by writing the Taylor expansion
of  $u$ at a given point $x\in\Omega$ and averaging it 
on a ball $B_\epsilon(x)$.  One obtains:
\begin{equation}\label{average1}
\fint_{B_\epsilon(x)} u = u(x) + \frac{\epsilon^2}{2(N+2)} \Delta u(x) +
o(\epsilon^2),
\end{equation}
which is a continuum version of (\ref{mean1})  when the second term in
the right hand side  vanishes. Consequently, the function $u$ must be {\em harmonic}, i.e. $\Delta u =0$.

\par\bigskip

\newpage

\noindent {\bf The $p$-Laplacian and its mean value property.}

To apply a similar reasoning to a nonlinear
problem, consider the homogeneous $p$-Laplacian:
\begin{equation}\label{plaplacian}
\Delta_p^{\textrm{H}} u =  \Delta u + (p-2 )\Delta_\infty u= |\nabla
u|^{2-p} \mbox{div}\big(|\nabla u|^{p-2} \nabla u\big),
\qquad 1<p<\infty, 
\end{equation}
where the  {\em infinity-Laplacian} is given by:  $\Delta_\infty
u = \langle \nabla^2u : \frac{\nabla u}{|\nabla u|} \otimes
\frac{\nabla u}{|\nabla u|}\rangle$.
Parallel to (\ref{average1}) one gets the expansion:
\begin{equation}\label{average2}
\frac{1}{2}\Big(\sup_{B_\epsilon(x)} u + \inf_{B_\epsilon(x)}
u\Big) = u(x) + \frac{\epsilon^2}{2}\Delta_\infty u(x) + o(\epsilon^2).
\end{equation}
Forming a linear combination of (\ref{average1}) and (\ref{average2}) with coefficients 
 $\alpha=
\frac{p-2}{p+N}$ and $\beta= \frac{2+N}{p+N}$, yields:
\begin{equation}\label{mean2}
\frac{\alpha}{2}\Big( \sup_{B_{\epsilon}(x)} u +
  \inf_{B_{\epsilon}(x)} u\Big) + \beta\fint_{B_\epsilon
  (x)} u = u(x) + \frac{\epsilon^2}{2(p+N)}\Delta_p^{\textrm{H}} u + 
o(\epsilon^2),
\end{equation}
and so the equation (\ref{mean2}) suggests that a {\em $p$-harmonic} function $u$,
i.e. a function satisfying $\Delta_p^{\textrm{H}} u  = 0$, may
be approximated by \textit{$p$-harmonious}  functions $u_{\epsilon}$,
defined by the mean value property:
\begin{equation}\label{mean3}
\frac{\alpha}{2}\Big( \sup_{B_{\epsilon}(x)} u_{\epsilon} +
  \inf_{B_{\epsilon}(x)} u_{\epsilon}\Big) + \beta\fint_{B_\epsilon
  (x)} u_{\epsilon} = u_{\epsilon}(x).
\end{equation}
As we shall see, the functions $u_\epsilon$ satisfying (\ref{mean3})
have a probabilistic interpretation as values of Tug-of-War  games
with noise. 

\par\bigskip 

\noindent {\bf A Tug-of-War game with noise for $\Delta_p^{\textrm{H}}$.}

\begin{figure}[h] 
{\includegraphics[trim = 25mm 63mm 25mm 70mm, clip, width=0.6\linewidth]{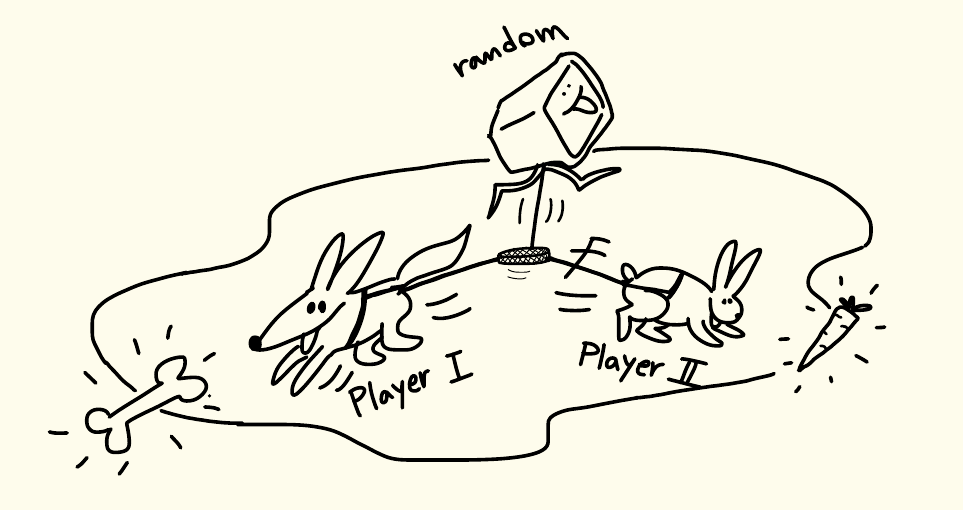}}
\caption{ Player I and Player II compete in a Tug-of-War with random
  noise \cite{Ka}.}
\label{figura1}
\end{figure}

\begin{figure}[h] 
\begin{picture}(150,130)
%\thinlines
\put(-50,60){\vector(1,1){50}}
\put(-50,60){\vector(1,0){50}}
\put(-50,60){\vector(1,-1){50}}
\put(-68,57){$\mathbf{x}_{n} $}
\put(-35,90){$\frac{\alpha}{2}$}
\put(-20,65){$\frac{\alpha}{2} $}
\put(-35,25){$\beta $}
\put(10,108){Player I chooses $\mathbf{x}_{n+1}$}
\put(10,58){Player II chooses $\mathbf{x}_{n+1}$}
\put(10,6){Random}
\put(130,84){\includegraphics[trim = 49mm 73mm 74mm 92mm, clip,
  width=0.13\linewidth]{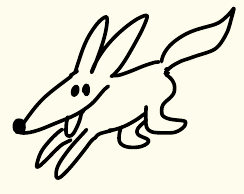}}
\put(130, 45){\includegraphics[trim = 87mm 78mm 45mm 93mm, clip,
  width=0.09\linewidth]{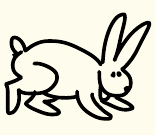}}
\put(130,-10){\includegraphics[trim = 73mm 91mm 59mm 74mm, clip,
  width=0.09\linewidth]{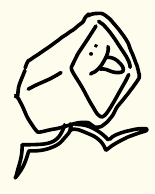}}
\end{picture}
\caption{ Player I, Player II and random noise with their probabilities \cite{Ka}.}
\label{figura2}
\end{figure}
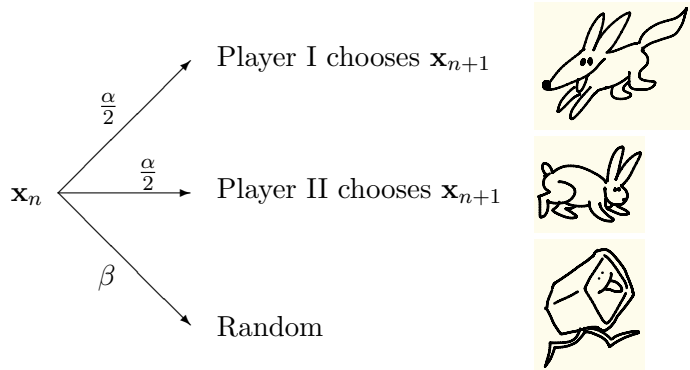

A Tug-of-War is a
two-person, zero-sum game, i.e. two players compete and the
gain of Player I equals the loss of Player II. Initially, a token is placed at a
point $x_0\in\Omega$. At each step of the process (the game) one of the
three actions takes place: (i) with probability
$\frac{\alpha}{2}$, Player I is allowed to play, and she moves the
token from its current position $x_n$ to her chosen position
$x_{n+1}\in B_\epsilon(x_n)$, (ii) with probability $\frac{\alpha}{2}$,
Player II moves the token to his chosen
position in $B_\epsilon(x_n)$, (iii) with probability $\beta=1-\alpha\in [0,1]$,
the token is moved randomly in the ball $B_\epsilon(x_n)$. The game
stops when the token leaves $\Omega$, whereas Player II pays to Player
I the amount equal to the value of a given boundary pay-off function $F$ at
the exit token position $x_\tau$ (see Figures \ref{figura1} and
\ref{figura2}).

%trim option's parameter order: left bottom right top

Players I and II play according to  {\em strategies} $\sigma_I$ and
$\sigma_{II}$ respectively, which are  (Borel measurable) functions assigning to each
finite history of the game ${\bf x}_n = (x_0,\ldots x_n)$ the next position $x_{n+1}$
in the $\epsilon$-neighborhood $\Omega_\epsilon$ of $\Omega$,
where the player will move the token once she/he is allowed to play.
These strategies determine  a probability
$\mathbb{P}^{x_0}_{\sigma_{I}, \sigma_{II}}$ on the space of all possible game runs in 
$(\Omega_\epsilon)^\infty$.  Since $\beta>0$, the game ends
$\mathbb{P}^{x_0}_{\sigma_{I}, \sigma_{II}}$-a.e., so that we can
define the stopping time $\tau(x_0, x_1,\ldots) = \inf\{n;~x_n\not\in\Omega\}$. 
The expected value of the game  is then given by $\mathbb{E}^{x_0}_{\sigma_{I}, \sigma_{II}}[F_\tau] =
\int_{\Omega^\infty } F(x_\tau) ~\mbox{d}\mathbb{P}^{x_0}_{\sigma_{I},
  \sigma_{II}}$. 
Consequently, the minimum gain $u_I$ that Player I can expect, and the
maximum loss $u_{II}$  of Player II  in his best case scenario, are:
\begin{equation}\label{gameval}
u_{I}(x_0) = \sup_{\sigma_{I}}\inf_{\sigma_{II}}
\mathbb{E}^{x_0}_{\sigma_{I}, \sigma_{II}}[F_\tau], \qquad 
u_{II}(x_0) = \inf_{\sigma_{II}}\sup_{\sigma_{I}}
\mathbb{E}^{x_0}_{\sigma_{I}, \sigma_{II}}[F_\tau].
\end{equation}
The following main results were achieved in \cite{PSSW} for
$p=\infty$,  and in \cite{MPR, LPS} for $p\in [2,\infty)$:

\begin{itemize}
\item[{}]  \textbf{Theorem A:} The two game values in (\ref{gameval}) coincide: $u_I=u_{II} $ and are equal
to the $p$-harmonious function $u_\epsilon$, which is the unique
solution to the mean value law in (\ref{mean3}) augmented by the boundary
data: $u_\epsilon = F$ on $\mathbb{R}^{n}\setminus \Omega$.

\medskip
 
\item[{}]  \textbf{Theorem B:}  As $\epsilon\to 0$, the game value $u_\epsilon$ converges uniformly in
$\Omega$ to a function $u\in\mathcal{C}(\Omega)$, which
is the unique viscosity solution to: $\Delta_p^{\textrm{H}} u =0$ in $\Omega$ and $u=F$ on $\partial\Omega$. 
 
\end{itemize}
In brief, one obtains exactly the same as in the case of the
discrete Brownian motion, whose value satisfied the averaging principle
(\ref{mean1}) and converged to a harmonic function. 
The equivalence
of the notions of  viscosity solution to $\Delta_p^{\textrm{H}} u=0$ and  weak solution to
$\mbox{div}\big(|\nabla u|^{p-2} \nabla u\big)=0$  has been proven in \cite{JLM}. 

\par\bigskip

\noindent{\bf The key martingale calculation.}

We now sketch the proof of
$u_{II}\leq u_\epsilon$; a symmetric argument yields $u_\epsilon\leq
u_{I}$, while $u_{I}\leq u_{II}$ is trivially  true. We take advantage of
the cancellation encoded in the mean value property (\ref{mean3}) by
showing that certain quantities related to $u_{\epsilon}$ are sub- and
super-martingales. 

Fix a small  $\eta>0$ and let $\bar\sigma_{II}$ be an \lq\lq almost
optimal\rq\rq  ~strategy for Player II, so that:
$$u_\epsilon(\bar\sigma_{II}(x_0,\ldots x_n)) \leq
\inf_{B_\epsilon(x_n)} u_\epsilon + \frac{\eta}{2^{n+1}}.$$
Player I plays according to an arbitrary strategy $\sigma_{I}$.
The key observation is that {\it the sequence of random variables 
$\{u_\epsilon(x_{n}) + \frac{\eta}{2^{n}}\mid (x_0, \ldots x_n)\}_{n\geq 1}$ is a
supermartingale}. 

%trim option's parameter order: left bottom right top

To see this, compute the conditional expectation relative to the history $\mathbf{x}_{n}=(x_{0}\ldots x_{n})$:
\begin{equation*}\label{mart}
\begin{split}
 \mathbb{E}^{x_0}_{\sigma_{I}, \bar\sigma_{II}}\big\{u_\epsilon(x_{n+1}) +
\frac{\eta}{2^{n+1}}\big\} (\mathbf{x}_{n} ) & =  \frac{\alpha}{2}
u_\epsilon(\sigma_I(\mathbf{x}_{n} ))+ \frac{\alpha}{2}
u_\epsilon(\bar\sigma_{II}(\mathbf{x}_{n} )) + \beta\fint_{B_\epsilon(x_n)}
u_\epsilon + \frac{\eta}{2^{n+1}}\\
& \leq \frac{\alpha}{2} \sup_{B_\epsilon(x_n)} u_\epsilon + \frac{\alpha}{2}
\left(\inf_{B_\epsilon(x_n)} u_\epsilon + \frac{\eta}{2^{n+1}}\right)
+ \beta\fint_{B_\epsilon(x_n)} u_\epsilon + \frac{\eta}{2^{n+1}}\\
& = u_\epsilon(x_n) + (\frac{\alpha}{2}+1) \frac{\eta}{2^{n+1}}\leq
 u_\epsilon(x_n) + \frac{\eta}{2^{n}},
\end{split}
\end{equation*}
where we first used the game's rules, then the
sub-optimality of $\bar\sigma_{II}$, and further the formula
(\ref{mean3}) for $u_\epsilon$.
Applying Doob's optimal stopping time theorem we get the desired comparison result:
\begin{equation*}\label{doob}
\begin{split}
u_{II}(x_0) & \leq \sup_{\sigma_{I}} 
\mathbb{E}^{x_0}_{\sigma_{I}, \bar\sigma_{II}}[F_\tau] = 
\sup_{\sigma_{I}} 
\mathbb{E}^{x_0}_{\sigma_{I}, \bar\sigma_{II}}[u(x_\tau)] 
\leq \sup_{\sigma_{I}} 
\mathbb{E}^{x_0}_{\sigma_{I}, \bar\sigma_{II}}[u(x_\tau) +
\frac{\eta}{2^\tau}] \\ &
\leq \sup_{\sigma_{I}} \mathbb{E}^{x_0}_{\sigma_{I}, \bar\sigma_{II}}[u(x_0) +
\frac{\eta}{2^0}] = u(x_0) + \eta,  \qquad \text{ for all } \eta>0.
\end{split}
\end{equation*}

\par\bigskip

\noindent{\bf Strategies and inequalities.}

We have  seen   how probability tools can be used 
%in a powerful way 
to study nonlinear PDEs, where the key technical ingredient was
assigning sutiable {\em strategies} yielding the desired {\em
  inequalities} for game values. Below we sketch two 
further examples of this powerful technique.

The proof of {\em uniform convergence} in Theorem B relies on a variant
of the Ascoli-Arzel\'a theorem valid for the discontinuous
functions $u_\epsilon$. The verification \cite{MPR}  of the appropriate
'equidiscontinuity' property requires estimating quantities
$|u_\epsilon(x_0) - u_\epsilon(y_0)|$, say for $x_0\in\Omega$, $y_0\in\partial\Omega$.
If $F$ is Lips\-chitz, this reduces to estimating $|x_\tau - y_0|$, and
the feasible strategy is that of Player II ``pulling towards $y_0$'', namely
shifting the token by $\epsilon$ along the segment connecting its current position
with $y_0$.

In \cite{LPSharnack}, the {\em local Harnack inequality} for $p$-harmonic
functions for $p>2$ is proven independent of the classical, yet technically challenging methods of
De Giorgi or Moser. This is done via a uniform estimate on the oscillations
of $u_\epsilon$. Let $x_0, y_0\in\Omega$ and let $z$ be equidistanced
from both points by a multiple of $\epsilon$. Define strategies
$\sigma^*_i$ in which Player $i$ cancells the earliest uncancelled
move of her/his opponent, and otherwise ``pulls towards $z$'' as
before, and let $\tau_{i}^*$ be the stopping time in which the game terminates
when either Player $i$ has played sufficiently many turns to place the
token at $z$ (modulo the random noise), or when the total amount of token's shifts
by her/his opponent and by the random noise, has passed an undesired
large treshold $r$.
Let now $\sigma_{I}, \sigma_{II}$ be two arbitrary strategies. By the
symmetry of this construction, the bulk ``nonlinear'' parts in the two quantities:
$\mathbb{E}_{\sigma_{I}, \sigma^*_{II}}^{x_0}[u_\epsilon(x_{\tau_{II}^*})]$ and
$\mathbb{E}_{\sigma^*_{I}, \sigma_{II}}^{y_0}[u_\epsilon(x_{\tau_{I}^*})]$,
corresponding to stopping the game due to the first condition, are
equal. The remaining ``linear'' part in: $|\mathbb{E}_{\sigma_{I},
  \sigma^*_{II}}^{x_0}[u_\epsilon(u_{\tau_{II}^*})] - \mathbb{E}_{\sigma^*_{I},
  \sigma_{II}}^{y_0}[u_\epsilon(u_{\tau_{I}^*})]|$ can then be bounded by $\frac{|x_0
  - y_0|}{r} \mbox{osc}(u_\epsilon, B_r(z))$, using a comparison with a
cylinder walk. This concludes the
proof, in view of: $|u_\epsilon(x_0) - u_\epsilon(y_0)|\leq \sup_{\sigma_I, \sigma_{II}}
|\mathbb{E}_{\sigma_{I},
  \sigma^*_{II}}^{x_0}[u_\epsilon(x_{\tau_{II}^*})] -
\mathbb{E}_{\sigma^*_{I},
  \sigma_{II}}^{y_0}[u_\epsilon(x_{\tau_{I}^*})]|$. 

\par\bigskip

\noindent{\bf Further results.}

Generalizations of Theorems A and B have been obtained in various
contexts. For $p=\infty$, only the notion of a metric space
is necessary to define the game, and indeed
\cite{PSSW} formulates its results for an arbitrary {\em length
  space} where the solutions to $\Delta_\infty u=0$ are understood as
{\it Absolutely Minimizing Lipschitz Extensions}.
When $p\in [2, \infty)$, the game uses  the notion of a metric and a measure, and it
is amenable to the the recent extension to Heisenberg groups  in \cite{FLM}.
For $p\in (1,\infty)$ one needs the additional notion of
perpendicularity \cite{PS}. 
We see that as $p\to 1$, the required complexity of structure increases.
In the case $p=1$ the game 
%is purely deterministic and it
is naturally related to the mean curvature flow \cite{KS} and functions of least
gradient.  

Other extensions include  the obstacle problems \cite{MRS}, finite difference
schemes \cite{AS}, equations with
right hand side $f\neq 0$, mixed boundary data \cite{APSS, CGR} and parabolic equations
\cite{KS2, MPR2}.


\begin{thebibliography}{9999}
\bibitem[Ka]{Ka} The cartoons in this paper were drawn by our colleague Kiumars
  Kaveh at the University of Pittsburgh.
\end{thebibliography}
\end{document}